\documentclass[leqno]{amsart}

\usepackage{amssymb}
\usepackage{amsfonts}
\usepackage{amsthm}

\theoremstyle{plain} 
\newtheorem{theorem}{Theorem}
\newtheorem{corollary}[theorem]{Corollary}
\newtheorem{lemma}[theorem]{Lemma}
\theoremstyle{definition} 

\theoremstyle{remark} 
\newtheorem{remark}[theorem]{Remark}

\newcommand{\R}{\ensuremath{\mathbb{R}}}
\newcommand{\Z}{\ensuremath{\mathbb{Z}}}
\newcommand{\T}{\ensuremath{\mathbb{T}}}

\DeclareMathOperator{\Log}{Log}
\DeclareMathOperator{\crd}{C_{rd}}
\DeclareMathOperator{\crdone}{C^1_{rd}}

\numberwithin{equation}{section}
\numberwithin{theorem}{section}

\small\normalsize

\begin{document}
\author[anderson]{Douglas R. Anderson}
\title[inhomogeneous second--order linear equations]{Alternative solutions of inhomogeneous second--order linear dynamic equations on time scales}
\address{Department of Mathematics and Computer Science, Concordia College, Moorhead, MN 56562 USA\\visiting the School of Mathematics and Statistics, The University of New South Wales, Sydney, NSW 2052, Australia}
\email{andersod@cord.edu}
\urladdr{http://www.cord.edu/faculty/andersod/bib.html}
\author[tisdell]{Christopher C. Tisdell}
\address{School of Mathematics and Statistics, The University of New South Wales, Sydney, NSW 2052, Australia}
\email{cct@unsw.edu.au}
\urladdr{http://web.maths.unsw.edu.au/~cct/}

\begin{abstract}
We exhibit an alternative method for solving inhomogeneous second--order linear ordinary dynamic equations on time scales, based on reduction of order rather than variation of parameters. Our form extends recent (and long-standing) analysis on $\R$ to a new form for difference equations, quantum equations, and arbitrary dynamic equations on time scales.
\end{abstract}

\maketitle

\thispagestyle{empty}


\section{Introduction to second-order ordinary dynamic equations}

A very common equation in mathematics, mathematical physics, and engineering is the inhomogeneous second--order linear ordinary differential equation
\begin{equation}\label{ode}
 y''(t)+p(t)y'(t)+q(t)y(t)=r(t), \quad t\in\R,
\end{equation}
its ordinary difference equation counterparts
\begin{equation}\label{odiffe}
 \Delta(\Delta y)(t)+p(t)\Delta  y(t)+q(t)y(t)=r(t), \quad t\in\Z, \quad \Delta y(t):=y(t+1)-y(t)
\end{equation}
or
\begin{equation}\label{odiffe2}
 y(t+2)+\alpha(t)y(t+1)+\beta(t)y(t)=r(t), \quad t\in\Z,
\end{equation}
the related ordinary quantum equation
\begin{equation}\label{oqe}
 D_h(D_h y)(t)+p(t)D_h y(t)+q(t)y(t)=r(t), \quad t\in h^\Z, \quad h>1, \quad D_h y(t):=\frac{y(ht)-y(t)}{ht-t},
\end{equation}
or the recently introduced ordinary dynamic equation on time scales given by
\begin{equation}\label{odyne}
 y^{\Delta\Delta}(t)+p(t)y^{\Delta}(t)+q(t)y(t)=r(t), \quad t\in\T^{\kappa^2},
\end{equation}
where the differential/shift operators in \eqref{ode}$-$\eqref{odyne} represent differentiation with respect to $t$ on the corresponding time scales, respectively. In the general setting represented by \eqref{odyne}, the functions $p$, $q$, and $r$ are real-valued right-dense continuous scalar functions of $t$ satisfying the regressivity condition
$$ 1+\mu(t)\left[-p(t)+\mu(t)q(t)\right]\neq 0. $$
Recall that on a time scale $\T$, namely any nonempty closed subset of the real line, the delta derivative is given by
$$ y^\Delta(t):=\lim_{s\rightarrow t}\frac{y^\sigma(t)-y(s)}{\sigma(t)-s}, \quad t\in\T^\kappa, $$
provided the limit exists, where $\sigma$ is the forward jump operator $\sigma(t):=\inf\{s\in\T: s>t\}$, and $y^\sigma=y\circ \sigma$. The graininess $\mu$ is simply $\mu(t)=\sigma(t)-t$. For more on time scales see Hilger \cite{hilger}.

Indeed, extensive analysis of \eqref{odyne} and its solution forms can be found in Bohner and Peterson \cite{bp1}, while the ordinary differential equation \eqref{ode} is studied by, for example, Blest \cite{blest}, Boyce and DiPrima \cite{bd}, Hille \cite{hille}, Ince \cite{ince}, Johnson, Busawon, and Barbot \cite{jbb}, Kelley and Peterson \cite{kp1}, whereas the difference equation \eqref{odiffe} appears in Agarwal \cite{ravi}, Elaydi \cite{elaydi}, and Kelley and Peterson \cite{kp2}.
A commonly used technique to solve \eqref{ode}$-$\eqref{odyne} is Lagrange's variation of parameters method. In this approach, a solution $y$ of \eqref{ode}$-$\eqref{odyne} takes the form $y=y_u+y_d$, where $y_u$ is the complementary solution of the corresponding homogeneous or undriven ($r\equiv 0$) form of \eqref{ode}$-$\eqref{odyne}, and $y_d$ is any particular solution of the inhomogeneous or driven equations \eqref{ode}$-$\eqref{odyne}. If $y_1$ and $y_2$ are two linearly independent solutions of the corresponding homogeneous equation, then it is well known that we may write $y=c_1y_1+c_2y_2+y_d$ for arbitrary constants $c_1$ and $c_2$. For example, using variation of parameters, the form of a particular solution for \eqref{odyne} is \cite[Theorem 3.73]{bp1}
\begin{equation}\label{bpvp}
 y_d(t) = y_2(t)\int_{a}^{t}\frac{y_1^\sigma(s)r(s)}{W^\sigma (y_1,y_2)(s)}\Delta s -  y_1(t)\int_{a}^{t}\frac{y_2^\sigma(s)r(s)}{W^\sigma (y_1,y_2)(s)}\Delta s, \quad t\in\T, 
\end{equation}
which reduces to 
$$ y_d(t) = y_2(t)\int_{a}^{t}\frac{y_1(s)r(s)}{W(y_1,y_2)(s)}ds - y_1(t)\int_{a}^{t}\frac{y_2(s)r(s)}{W(y_1,y_2)(s)}ds, \quad t\in\R $$
for \eqref{ode}, and to
$$ y_d(t) = y_2(t)\sum_{s=a}^{t-1}\frac{y_1(s+1)r(s)}{W(y_1,y_2)(s+1)} - y_1(t)\sum_{s=a}^{t-1}\frac{y_2(s+1)r(s)}{W(y_1,y_2)(s+1)}, \quad t\in\Z $$
for \eqref{odiffe}, where in each case $W(y_1,y_2)$ is the Wronskian of $y_1$ and $y_2$, defined appropriately for each time scale. Notice that the integrals/summations always involve both $y_1$ and $y_2$. Recently, Johnson, Busawon, and Barbot \cite{jbb} derived two alternative solution forms for $y_d$ for the ordinary differential equation \eqref{ode}, namely
\begin{equation}\label{jbbform}
 y_d(t)=y_i(t)\int \frac{e^{-\int p(t)dt}\int r(t)y_i(t)e^{\int p(t)dt}dt}{y_i^2(t)}dt, \quad t\in\R, \quad i=1,2.
\end{equation}
As the authors point out in \cite{jbb}, either $y_1$ or $y_2$ may be chosen in \eqref{jbbform}, depending on which one yields an easier integral to compute. 

Unfortunately neither the technique nor the results in \cite{jbb} are new, as they are discussed in Blest \cite{blest} and the undergraduate textbook by Boyce and DiPrima \cite[Exercise 28, p185]{bd}, to cite just two examples. In fact, Yosida \cite[p.29]{yosida} calls the technique D'Alembert's reduction of order for ordinary differential equations, while according to Jahnke \cite[p.332]{jahnke}, the method of reduction of order dates back at least as far as Euler (1750). In addition, Clairaut, Lagrange and Laplace have been involved with the method and thus it may be impossible to attribute to one mathematician \cite[p.332]{jahnke}. In any case, the method is at least 250 years old.

Despite the non-novelty of \eqref{jbbform} for the ode \eqref{ode}, such a result would be new on general time scales for \eqref{odyne}. Consequently one of our goals in this paper is to generalize and extend \eqref{jbbform} to an alternative form for particular solutions to the ordinary dynamic equation on time scales \eqref{odyne}, which would then nevertheless result in new alternative forms for the difference equations \eqref{odiffe} or \eqref{odiffe2}, and the quantum equation \eqref{oqe} as simple corollaries, as well as contain \eqref{jbbform} for \eqref{ode}. This development is different from that given by Bohner and Peterson \cite[Chapter 3]{bp1}.

An additional goal of this paper is to present some simple results that guarantee the nonmultiplicity of solutions to initial value problems associated with \eqref{odyne}.  Our approach is based on simple inequalities and does not require a knowledge of matrix theory, in contrast to \cite[Theorem 3.1; Corollary 5.90]{bp1}.

In Section \ref{sec2} we address the question of nonmultiplicity of solutions, while in Section \ref{sec3} we develop a method for solving inhomogeneous second--order linear ordinary dynamic equations on time scales, based on reduction of order rather than variation of parameters.  Section \ref{sec4} contains some special cases that illustrate our results.


\section{nonmultiplicity} \label{sec2}

In this section we consider the notion of nonmultiplicity of solutions to the linear initial value problem \eqref{odyne} with initial conditions 
\begin{equation}\label{ics}
 y(t_0)=A, \qquad y^\Delta(t_0)=B,
\end{equation}
where $A,B\in\R$ and $t_0\in\T$.  Let $I \subseteq \R$ and $t_0 \in I_\T := I \cap \T$.  
 By ``nonmultiplicity of solutions'' we mean that our theorems will present conditions under which \eqref{odyne}, \eqref{ics} will have, at most, one solution $y = y(t)$ for $t\ge t_0$, $t\in I_\T$. 
Such information is highly valuable, for example, when constructing explicit solutions to problems as we can determine when the constructed solution (or unique linear combination of solutions) will be the only solution to the problem at hand.  Our techniques follow those of Coddington \cite[Ch.2, Sec.3]{coddington}.

In the following, we will make use of the time-scale exponential function defined in terms of right-dense continuous functions $p$ satisfying the regressivity condition $1+\mu(t)p(t)\neq 0$ for all $t\in\T^\kappa$. Given such a $p$, the delta exponential function \cite[Theorem 2.30]{bp1} is given by 
$$ e_p(t,a) = \begin{cases} \exp\left(\displaystyle\int_a^t p(\tau)\Delta\tau\right) &: \mu(\tau) = 0, \\
  \exp\left(\displaystyle\int_a^t\frac{1}{\mu(\tau)}\Log(1+p(\tau)\mu(\tau))\Delta\tau\right) &: \mu(\tau)\ne 0, \end{cases} $$
where $\Log$ is the principal logarithm. It follows that $e_p(t,a)$ is the unique solution to the initial value problem $\phi^\Delta(t)=p(t)\phi(t)$, $\phi(a)=1$ on $\T$. We will denote $1/e_p(t,a)$ by $e_{\ominus p}(t,a)$. 

We will require the following lemma from \cite[Theorem 6.1, p.255]{bp1}. 


\begin{lemma}\label{cctlemma1}
Let $\ell\in\mathcal{R}^+$ and $v\in\crdone(\T)$, and let $t_0\in\T$. If $v^\Delta(t) \le \ell(t)v(t)$ for all $t\ge t_0$, then $v(t) \le v(t_0)e_{\ell}(t,t_0)$ for all $t\ge t_0$.
\end{lemma}

Our initial analysis will concern the case of \eqref{odyne} with constant coefficients, namely
\begin{equation}\label{cct3}
 y^{\Delta\Delta}(t)+py^{\Delta}(t)+qy(t)=r(t), \quad y(t_0)=A, \quad y^\Delta(t_0)=B,
\end{equation}
where $p,q\in\R$ are constants. The following result involves the homogeneous form of \eqref{cct3}, that is
\begin{equation}\label{cct4}
 y^{\Delta\Delta}(t)+py^{\Delta}(t)+qy(t)=0,
\end{equation}
and provides an estimate on the growth rate of solutions to \eqref{cct4}. For this estimate, we define
$$ \|y(t)\|_2:=\left((y(t))^2+\left(y^\Delta(t)\right)^2\right)^{1/2}, \quad t\in I_\T, \quad t\ge t_0. $$


\begin{theorem}\label{ccthm1}
Let $k:=1+|p|+|q|$. If $y$ is any solution to \eqref{cct4} on $I_\T$ then
\begin{equation}\label{cct5}
 \|y(t)\|_2 \le \|y(t_0)\|_2e_k(t,t_0), \quad t\in I_\T, \quad t\ge t_0.
\end{equation}
\end{theorem}

\begin{proof}
Let $u(t)=\|y(t)\|_2^2$, where $y$ is a solution to \eqref{cct4}. For all $t\in I_{\T}^\kappa$ we have
\begin{eqnarray*}
 u^\Delta(t) &=& \left(y(t)+y^\sigma(t)\right)y^\Delta(t) + \left(y^\Delta(t)+y^{\Delta\sigma}(t)\right)y^{\Delta\Delta}(t) \\
 &=& \left(2y(t)+\mu(t)y^\Delta(t)\right)y^\Delta(t) + \left(2y^\Delta(t)+\mu(t)y^{\Delta\Delta}(t)\right)y^{\Delta\Delta}(t) \\
 &=& 2y(t)y^\Delta(t) + \mu(t)(y^\Delta(t))^2 + 2y^\Delta(t)y^{\Delta\Delta}(t) + \mu(t)(y^{\Delta\Delta}(t))^2.
\end{eqnarray*}
Now, apply Young's inequality $2ab \le a^2 + b^2$ to the first term and replace $y^{\Delta\Delta}$ with $-py^\Delta-qy$ to obtain
\begin{eqnarray*}
 u^\Delta(t) &\le& (y(t))^2 + (y^\Delta(t))^2 + \mu(t)(y^\Delta(t))^2 + 2y^\Delta(t)[-py^\Delta(t)-qy(t)] + \mu(t)[-py^\Delta(t)-qy(t)]^2 \\
 &\le& (1+2|p|+|q|)\left((y(t))^2 + (y^\Delta(t))^2\right) + \mu(t)\left[1+p^2+|p||q|+q^2\right]\left((y(t))^2 + (y^\Delta(t))^2\right) \\
 &\le& 2(1+|p|+|q|)\left((y(t))^2 + (y^\Delta(t))^2\right) + \mu(t)\left[1+2(|p|+|q|)+p^2+2|p||q|+q^2\right]\left((y(t))^2 + (y^\Delta(t))^2\right) \\
 & = & (k\oplus k)u(t).
\end{eqnarray*}
Above, $\oplus$ is known as the ``circle plus'' operator, $z \oplus w := z + w + \mu zw$, see \cite[p.54]{bp1}.

Thus, the conditions of Lemma \ref{cctlemma1} hold with $v=u$ and $\ell=k\oplus k$. Consequently, 
$$ u(t) \le u(t_0)e_{k\oplus k}(t,t_0) = u(t_0) (e_k(t,t_0))^2, \quad t\ge t_0, \quad t\in I_\T, $$
and therefore
$$ \|y(t)\|_2 \le \|y(t_0)\|_2 e_{k}(t,t_0), \quad t\ge t_0, \quad t\in I_\T. $$
\end{proof}

Theorem \ref{ccthm1} now leads to the following nonmultiplicity result for solutions to \eqref{cct3}, \eqref{ics}.


\begin{theorem}\label{ccthm2}
Let $r\in\crd(I_{\T}^{\kappa^2};\R)$, and $t_0\in I_{\T}$. The dynamic initial value problem \eqref{cct3}, \eqref{ics} has, at most, one solution $y=y(t)$ for $t\ge t_0$ with $t\in I_{\T}$.
\end{theorem}

\begin{proof}
Assume there are two solutions $x$ and $y$. Let $z(t):=x(t)-y(t)$ for $t\in I_\T$, and note that $z$ must satisfy the homogeneous equation \eqref{cct4} together with the homogeneous initial conditions $z(t_0)=0$, $z^\Delta(t_0)=0$. Now, by Theorem \ref{ccthm1} we have
$$ \|z(t)\|_2 \le \|z(t_0)\|_2 e_k(t,t_0)=0, \quad t\ge t_0, \quad t\in I_{\T}, $$
and thus $z(t)=0$ for all $t\ge t_0$, $t\in I_{\T}$. It follows that $x=y$.
\end{proof}

The following result is an extension of Theorem \ref{ccthm1} and concerns estimates on solutions to the homogeneous form of \eqref{odyne} with variable coefficients, namely \begin{equation}\label{homot}
 y^{\Delta\Delta}(t)+p(t)y^{\Delta}(t)+q(t)y(t)=0, \quad t\in\T^{\kappa^2}.
\end{equation}


\begin{theorem}\label{ccthm3}
Let $t\in I_{\T}$. Let $p_1,q_1\in\R$ be nonnegative constants such that
$$ |p(t)|\le p_1, \quad |q(t)|\le q_1, \quad t\in I_{\T}^{\kappa^2}, \quad t\ge t_0, $$
and let $k_1:=1+p_1+q_1$. If $y$ is any solution to \eqref{homot} on $I_{\T}$, then 
$$ \|y(t)\|_2 \le \|y(t_0)\|_2 e_{k_1}(t,t_0), \quad t\ge t_0, \quad t\in I_{\T}. $$
\end{theorem}

\begin{proof}
Our proof follows similar lines as that of the proof of Theorem \ref{ccthm1}, and thus we just sketch the details. Letting $u(t)=\|y(t)\|_2^2$, we obtain
\begin{eqnarray*}
 u^\Delta(t) &\le& (y(t))^2 + (y^\Delta(t))^2 + \mu(t)(y^\Delta(t))^2 + 2y^\Delta(t)[-p(t)y^\Delta(t)-q(t)y(t)] + \mu(t)[-p(t)y^\Delta(t)-q(t)y(t)]^2 \\
 &\le& 2(1+p_1+q_1)u(t) + \mu(t)\left[1+2(p_1+q_1)+p_1^2+2p_1q_1+q_1^2\right]u(t) \\
 &\le& (k_1\oplus k_1)u(t).
\end{eqnarray*}
Thus, applying Lemma \ref{cctlemma1} we obtain 
$$ u(t) \le u(t_0) (e_{k_1}(t,t_0))^2, \quad t\ge t_0, \quad t\in I_\T, $$
and therefore
$$ \|y(t)\|_2 \le \|y(t_0)\|_2 e_{k_1}(t,t_0), \quad t\ge t_0, \quad t\in I_\T. $$
\end{proof}

Theorem \ref{ccthm3} now leads to the following nonmultiplicity  result for solutions to \eqref{odyne}, \eqref{ics}.


\begin{theorem}\label{ccthm4}
Let $t\in I_{\T}$. The dynamic equation \eqref{odyne} with initial conditions \eqref{ics} has at most one solution $y=y(t)$ for $t\ge t_0$ with $t\in I_{\T}$.
\end{theorem}

\begin{proof}
Assume there are two solutions $x$ and $y$, and let $z(t)=x(t)-y(t)$ for $t\ge t_0$ with $t\in I_{\T}$. Note that $z$ satisfies \eqref{homot} and the initial homogeneous initial conditions $z(t_0)=0$, $z^\Delta(t_0)=0$. Since $I_{\T}$ may be unbounded, the coefficient functions $p$ and $q$ may not be bounded on $I_{\T}$, and so Theorem \ref{ccthm3} may not be directly applied to $z$. We let $t$ be any point in $I_{\T}$ such that $t>t_0$, and let $J_{\T}$ be any closed, bounded interval of $I_{\T}$ such that $J_{\T}$ has $t_0$ as a left endpoint and $J_{\T}$ contains $t$. On $J_{\T}$, $p$ and $q$ are both bounded, say by $p_1$ and $q_1$, respectively. We can now apply Theorem \ref{ccthm3} to $z$ on $J_{\T}$, and so $z(t)=0$ for all $t\in J_{\T}$. Now since $t$ was chosen to be any point in $I_{\T}$ with $t>t_0$, we have shown that $x(t)=y(t)$ for all $t\in I_{\T}$ with $t\ge t_0$.
\end{proof}


\begin{remark}
The quest for nonmultiplicity of solutions on intervals to the left of $t_0$ is a more delicate affair.  The results of this section may be extended to include nonmultiplicity of solution for $t\le t_0$ by adapting the proofs and obtaining inequalities like
$$u^\Delta(t) \ge -2k u(t), \quad \mbox{for all} \quad t \le t_0, \ t \in I_\T.$$
However, there is a price to pay -- the graininess function of the time scale would need to be bounded above.  This is due to regressivity coming into play. 
\end{remark}


\section{alternative solution forms for the general inhomogeneous dynamic equation} \label{sec3}

In this section we state and prove the main result, namely a new form for a solution of \eqref{odyne} in the spirit of \eqref{jbbform}. As mentioned previously, this will then give us a new form for the difference equation \eqref{odiffe} and quantum equation \eqref{oqe} as well. 


\begin{theorem}\label{thmT}
Let $\T$ be a time scale and let $p$, $q$ and $r$ be real--valued right--dense continuous scalar functions of $t\in\T$ with $p$ and $q$ satisfying the regressivity condition
\begin{equation}\label{regT}
 1+\mu(t)\left[-p(t)+\mu(t)q(t)\right]\neq 0, \quad t\in\T^\kappa.
\end{equation}
For all $t \in \T$, let $y_1$ and $y_2$ satisfy
\begin{equation}\label{homoslomo}
 y_i^{\Delta\Delta}(t)+p(t)y_i^\Delta(t)+q(t)y_i(t)=0, \quad i=1,2, 
\end{equation}
and $W(t):=y_1(t)y_2^\Delta(t)-y_2(t)y_1^\Delta(t)\neq 0$.
The general solution of the linear inhomogeneous second--order ordinary dynamic equation
$$ y^{\Delta\Delta}(t)+p(t)y^{\Delta}(t)+q(t)y(t)=r(t), \quad t\in\T^{\kappa^2}, $$
is
\begin{equation}\label{odynform1}
 y(t)=c_1y_1(t)+c_2y_2(t)+y_1(t)\int \frac{e_{(-p+\mu q)}(t,a)\int r(t)y_1^\sigma(t)e_{\ominus(-p+\mu q)}(\sigma(t),a) \Delta t} {y_1(t)y_1^\sigma(t)}\Delta t
\end{equation}
or
\begin{equation}\label{odynform2}
 y(t)=c_1y_1(t)+c_2y_2(t)+y_2(t)\int \frac{e_{(-p+\mu q)}(t,a)\int r(t)y_2^\sigma(t)e_{\ominus(-p+\mu q)}(\sigma(t),a) \Delta t} {y_2(t)y_2^\sigma(t)}\Delta t,
\end{equation}
where $c_1$ and $c_2$ are arbitrary constants and $e_{(-p+\mu q)}(\cdot,a)$ is the time--scale exponential function.
\end{theorem}

\begin{proof}
Clearly \eqref{odynform1} and \eqref{odynform2} are of the expected form $y=c_1y_1+c_2y_2+y_d$. One could easily verify that a function of the form
\begin{equation}\label{ydform}
 y_d(t)=y_i(t)\int \frac{e_{(-p+\mu q)}(t,a)\int r(t)y_i^\sigma(t)e_{\ominus(-p+\mu q)}(\sigma(t),a) \Delta t} {y_i(t)y_i^\sigma(t)}\Delta t
\end{equation}
is a particular solution of \eqref{odyne} using the time scale calculus, but that would not give any insight into where \eqref{ydform} comes from. Thus, to derive \eqref{ydform}, assume we have a particular solution to the inhomogeneous equation \eqref{odyne} of the form $y_d(t)=y_i(t)v(t)$, where $y_i$ solves the homogeneous equation \eqref{homoslomo} and $v$ is a function to be determined. Then, using the product rule 
$(fg)^\Delta = gf^\Delta + f^\sigma g^\Delta$, we have
\begin{eqnarray*}
 y_d^\Delta(t) &=& v(t)y_i^\Delta(t)+y_i^\sigma(t)v^\Delta(t) \\
  &=& v(t)y_i^\Delta(t)+y_i(t)y_i^\sigma(t)v^\Delta(t)/y_i(t),
\end{eqnarray*}
and using the product rule again with the quotient rule $\left(\frac{f}{g}\right)^\Delta=\frac{gf^\Delta-fg^\Delta}{gg^\sigma}$ we see that
\begin{eqnarray*}
 y_d^{\Delta\Delta}(t) &=& v(t)y_i^{\Delta\Delta}(t)+y_i^{\Delta\sigma}(t)v^\Delta(t)+\frac{y_i(t)\left(y_i(t)y_i^\sigma(t)v^\Delta(t)\right)^\Delta-y_i(t)y_i^\sigma(t)v^\Delta(t)y_i^\Delta(t)}{y_i(t)y_i^\sigma(t)} \\
  &=& v(t)y_i^{\Delta\Delta}(t)+y_i^{\Delta\sigma}(t)v^\Delta(t)+\left(y_i(t)y_i^\sigma(t)v^\Delta(t)\right)^\Delta/y_i^\sigma(t)-v^\Delta(t)y_i^\Delta(t).
\end{eqnarray*}
Since we are assuming $y_d$ is a particular solution of the inhomogeneous equation \eqref{odyne}, we must have
\begin{eqnarray*}
 r(t) &=& y_d^{\Delta\Delta}(t)+p(t)y_d^{\Delta}(t)+q(t)y_d(t) \\
  &=& v(t)\left(y_i^{\Delta\Delta}(t)+p(t)y_i^{\Delta}(t)+q(t)y_i(t)\right)+p(t)y_i(t)y_i^\sigma(t)v^\Delta(t)/y_i(t) \\
  & & +y_i^{\Delta\sigma}(t)v^\Delta(t)+\left[y_i(t)y_i^\sigma(t)v^\Delta(t)\right]^\Delta/y_i^\sigma(t)-v^\Delta(t)y_i^\Delta(t).
\end{eqnarray*}
Now $y_i$ is a solution of the homogeneous equation \eqref{homoslomo}, so after simplifying we multiply by $y_i^\sigma$ to get
$$ r(t)y_i^\sigma(t)=y_i^\sigma(t)p(t)\left[y_i(t)y_i^\sigma(t)v^\Delta(t)\right]/y_i(t)+\left[y_i(t)y_i^\sigma(t)v^\Delta(t)\right]^\Delta
+y_i^\sigma(t)v^\Delta(t)\left(y_i^{\Delta\sigma}(t)-y_i^\Delta(t)\right). $$
Make the substitution $u=y_iy_i^\sigma v^\Delta$, then use the simple formula $f^\sigma-f=\mu f^\Delta$ to get
\begin{eqnarray*}
 r(t)y_i^\sigma(t) &=& y_i^\sigma(t)p(t)u(t)/y_i(t) + u^\Delta(t)+y_i^\sigma(t)v^\Delta(t)\mu(t)y_i^{\Delta\Delta}(t) \\
 &=& u^\Delta(t)+\frac{p(t)y_i^\sigma(t)+\mu(t)y_i^{\Delta\Delta}(t)}{y_i(t)}u(t).
\end{eqnarray*}
Using the simple formula $f^\sigma-f=\mu f^\Delta$ again and rearranging, we see that
$$ u^\Delta(t) + \left(p(t)-\mu(t)q(t)\right)u(t) = r(t)y_i^\sigma(t). $$
Focusing on the coefficient of $u$, we note that 
\begin{eqnarray*}
 p(t)-\mu(t) q(t) &=& -(-p(t)+\mu(t) q(t)) \\
 &=& -\ominus\left(\ominus(-p(t)+\mu(t) q(t))\right) \\
 &=& \frac{\ominus(-p(t)+\mu(t) q(t))e_{\ominus(-p+\mu q)}(t,a)}{\left[1+\mu(t)\left(\ominus(-p(t)+\mu(t) q(t))\right)\right]e_{\ominus(-p+\mu q)}(t,a)}\\
 &=& \frac{e^\Delta_{\ominus(-p+\mu q)}(t,a)}{e^\sigma_{\ominus(-p+\mu q)}(t,a)}.
\end{eqnarray*}
Consequently we have that
$$ \left(e_{\ominus(-p+\mu q)}(t,a)u(t)\right)^\Delta = r(t)y_i^\sigma(t)e^\sigma_{\ominus(-p+\mu q)}(t,a). $$
Solving while recalling that $u=y_iy_i^\sigma v^\Delta$, we arrive at
$$ v(t)=\int \frac{e_{(-p+\mu q)}(t,a)\int r(t)y_i^\sigma(t)e_{\ominus(-p+\mu q)}(\sigma(t),a) \Delta t} {y_i(t)y_i^\sigma(t)}\Delta t, $$
so that via $y_d=y_iv$ we obtain \eqref{ydform}.
\end{proof}

\begin{remark}
The regressivity assumption in \eqref{regT} is not at all unusual, as it is automatic in the case $\T=\R$ since $\mu\equiv 0$, and it is assumed in the variation of parameters theorem on general time scales; see \cite[Definition 3.3]{bp1}. The following theorem is a simple corollary of Theorem \ref{thmT}.
\end{remark}

\begin{theorem}[Reduction of Order]\label{roo}
Assume $p$ and $q$ are real--valued right--dense continuous scalar functions of $t\in\T$ satisfying the regressivity condition
\begin{equation}\label{regTpos}
 1+\mu(t)\left[-p(t)+\mu(t)q(t)\right]\neq 0, \quad t\in\T^{\kappa}.
\end{equation}
If $y_1$ is a solution of the linear homogeneous second--order ordinary dynamic equation \eqref{homot}, then
\begin{equation}\label{rooy2}
 y_2(t)=y_1(t)\int \frac{e_{(-p+\mu q)}(t,a)}{y_1(t)y_1^\sigma(t)}\Delta t
\end{equation}
is a second linearly independent solution of \eqref{homot}. Similarly, if $y_2$ is a solution of \eqref{homot}, then
\begin{equation}\label{rooy1}
 y_1(t)=y_2(t)\int \frac{e_{(-p+\mu q)}(t,a)}{y_2(t)y_2^\sigma(t)}\Delta t
\end{equation}
is a second linearly independent solution of \eqref{homot}. 
\end{theorem}

\begin{proof}
We will prove \eqref{rooy2}, since the proof of \eqref{rooy1} is similar. Thus, assume $y_1$ is a solution of \eqref{homot}. Since $r(t)\equiv 0$ in this case, and general antiderivatives are used in Theorem \ref{thmT}, we may choose the constant of integration in \eqref{ydform} in such a way that
$$ \int r(t)y_1^\sigma(t)e_{\ominus(-p+\mu q)}(\sigma(t),a) \Delta t=1, $$
so that \eqref{ydform} becomes
$$ y_d(t)=y_1(t)\int \frac{e_{(-p+\mu q)}(t,a)}{y_1(t)y_1^\sigma(t)}\Delta t=y_2(t), $$
in other words a particular solution of \eqref{homot}. One could also verify \eqref{rooy2} directly using the delta derivative rules. To show linear independence, we calculate the Wronskian of $y_1$ and $y_2$, namely
\begin{equation}\label{Wron12}
 W(y_1,y_2)(t) = y_1(t) y_2^\Delta(t) - y_1^\Delta(t) y_2(t) = e_{(-p+\mu q)}(t,a) \neq 0 
\end{equation}
for all $t\in\T$ by the regressivity assumption \eqref{regTpos}; see \cite[Theorem 2.48]{bp1}. 
\end{proof}

\begin{remark}
It is also possible to show the equivalence of the result shown in \eqref{ydform} with that due to the variation of parameters formula in equation \eqref{bpvp}. Let $y_1$ be a solution of the linear homogeneous equation \eqref{homot}. Assuming \eqref{regTpos}, and appropriating the Wronskian of $y_1$ and the form of $y_2$ given in \eqref{rooy2} as calculated in \eqref{Wron12}, the solution given in \eqref{bpvp} is
\begin{eqnarray}
 y_d(t) &=& y_2(t)\int_{a}^{t}\frac{y_1^\sigma(s)r(s)}{e_{(-p+\mu q)}(\sigma(s),a)}\Delta s - y_1(t)\int_{a}^{t}\frac{y_2^\sigma(s)r(s)}{e_{(-p+\mu q)}(\sigma(s),a)}\Delta s \nonumber  \\
 &=& \left(y_1(t)\int_{a}^{t} \frac{e_{(-p+\mu q)}(s,a)}{y_1(s)y_1^\sigma(s)}\Delta s\right)\int_{a}^{t} r(s)y_1^\sigma(s)e_{\ominus(-p+\mu q)}(\sigma(s),a)\Delta s \nonumber \\
 & & -y_1(t)\int_{a}^{t}r(s)y_1^\sigma(s)e_{\ominus(-p+\mu q)}(\sigma(s),a)\int_{a}^{\sigma(s)}
 \frac{e_{(-p+\mu q)}(\xi,a)}{y_1(\xi)y_1^\sigma(\xi)}\Delta \xi\Delta s.\label{ydequivvp}
\end{eqnarray}
If we use the integration by parts formula $\int fg^\Delta\Delta t=f(t)g(t)-\int f^\Delta g^\sigma\Delta t$, see \cite[Theorem 1.77(vi)]{bp1}, on the integral in \eqref{ydform}, where we have taken
$$ f(t)=\int_{a}^{t}r(s)y_1^\sigma(s)e_{\ominus(-p+\mu q)}(\sigma(s),a)\Delta s \quad\text{and}\quad g^\Delta(t)=\frac{e_{(-p+\mu q)}(t,a)}{y_1(t)y_1^\sigma(t)}, $$ 
we get \eqref{ydequivvp}.
\end{remark}

The next corollary concerns another possible second--order linear dynamic equation discussed by Bohner and Peterson \cite[(3.6)]{bp1}.


\begin{corollary}\label{corT}
Let $\alpha$, $\beta$ and $r$ be real--valued right--dense continuous scalar functions on $\T$ with $\alpha$ satisfying the regressivity condition $1+\mu(t)\alpha(t)\neq 0$ for $t\in\T^{\kappa}$. For all $t \in \T$, let
 $y_1$ and $y_2$ satisfy 
\begin{equation}\label{corhomoslomo}
 y_i^{\Delta\Delta}(t)+\alpha(t)y_i^{\Delta\sigma}(t)+\beta(t)y_i^\sigma(t)=0, \quad i=1,2, 
\end{equation}
with $W(t):=y_1(t)y_2^\Delta(t)-y_2(t)y_1^\Delta(t)\neq 0$. The general solution of the linear inhomogeneous second--order ordinary dynamic equation
\begin{equation}\label{bp3.6}
 y^{\Delta\Delta}(t)+\alpha(t)y^{\Delta\sigma}(t)+\beta(t)y^\sigma(t)=r(t), \quad t\in\T^{\kappa^2}, 
\end{equation}
is
\begin{equation}\label{corodynform1}
 y(t)=c_1y_1(t)+c_2y_2(t)+y_1(t)\int\frac{e_{\ominus\alpha}(t,a)\int r(t)y_1^\sigma(t)e_{\alpha}(\sigma(t),a) \Delta t}{y_1(t)y_1^\sigma(t)}\Delta t
\end{equation}
or
\begin{equation}\label{corodynform2}
 y(t)=c_1y_1(t)+c_2y_2(t)+y_2(t)\int\frac{e_{\ominus\alpha}(t,a)\int r(t)y_2^\sigma(t)e_{\alpha}(\sigma(t),a) \Delta t}{y_2(t)y_2^\sigma(t)}\Delta t,
\end{equation}
where $c_1$ and $c_2$ are arbitrary constants, and where $e_{\alpha}(\cdot,a)$ is the time-scale exponential function.
\end{corollary}

\begin{proof}
Rewriting \eqref{bp3.6} using the simple formula $f^\sigma=f+\mu f^\Delta$, we see that we arrive at an equation of the form \eqref{odyne}, where
$$ p(t)=\frac{\alpha(t)+\mu(t)\beta(t)}{1+\mu(t)\alpha(t)}  \quad\text{and}\quad  q(t)=\frac{\beta(t)}{1+\mu(t)\alpha(t)}. $$
Then we see that the term $-p+\mu q$ in Theorem \ref{thmT} above is given by
$$ -p(t)+\mu(t) q(t) = \frac{-\alpha(t)}{1+\mu(t)\alpha(t)}=\ominus\alpha(t), $$
and this corollary follows.
\end{proof}


\section{second--order linear ordinary difference equations} \label{sec4}

In this section we recount the results of the previous section when $\T=\Z$, that is to say for second--order linear ordinary difference equations.
The following corollaries to Theorem \ref{thmT} are obtained by simply taking $\T=\Z$ in Theorem \ref{thmT} above. The first uses the forward difference operator form of the second--order linear equation \eqref{odiffe}, and the second uses the shift form \eqref{odiffe2}.


\begin{theorem}\label{resultz}
Let $p$, $q$ and $r$ be real--valued scalar functions of $t\in\Z$ with $p$ and $q$ satisfying the regressivity condition
\begin{equation}
 1-p(t)+q(t)\neq 0,  \quad t\in\Z.
\end{equation}
Let $y_1$ and $y_2$ satisfy 
$$ \Delta^2y_i(t)+p(t)\Delta  y_i(t)+q(t)y_i(t)=0, \quad i=1,2, $$
and $W(t):=y_1(t)\Delta  y_2(t)-y_2(t)\Delta  y_1(t)\neq 0$.
The general solution of the linear inhomogeneous second--order ordinary difference equation
$$ \Delta^2y(t) + p(t)\Delta y(t) + q(t)y(t) = r(t), \quad t\in\Z, $$
is given by
\begin{equation}\label{odiffform1} 
 y(t) = c_1 y_1(t) + c_2 y_2(t) + y_1(t) \sum \frac{\prod_{j=a}^{t-1} \left(1-p(j)+q(j)\right) \sum \frac{r(t)y_1(t+1)} {\prod_{j=a}^{t}\left(1-p(j)+q(j)\right)}} {y_1(t)y_1(t+1)}
\end{equation}
or
\begin{equation}\label{odiffform2}
y(t)=c_1y_1(t)+c_2y_2(t)+y_2(t)\sum\frac{\prod_{j=a}^{t-1}\left(1-p(j)+q(j)\right)\sum\frac{r(t)y_2(t+1)}{\prod_{j=a}^{t}\left(1-p(j)+q(j)\right)}}{y_2(t)y_2(t+1)},
\end{equation}
and where $c_1$ and $c_2$ are arbitrary constants.
\end{theorem}


\begin{theorem}\label{resultz2}
Let $\alpha$, $\beta$ and $r$ be real--valued scalar functions of $t\in\Z$ with regressivity condition $\beta(t)\neq 0$ for $t\in\Z$.
Let $y_1$ and $y_2$ satisfy 
$$ y_i(t+2)+\alpha(t)y_i(t+1)+\beta(t)y_i(t)=0, \quad i=1,2, $$
with $W(t):=y_1(t)y_2(t+1)-y_2(t)y_1(t+1)\neq 0$.
The general solution of the linear inhomogeneous second--order ordinary difference equation
\begin{equation}\label{elaydiform}
 y(t+2) + \alpha(t)y(t+1) + \beta(t)y(t) = r(t), \quad t\in\Z, 
\end{equation}
is given by
\begin{equation}\label{odiffform12} 
 y(t) = c_1 y_1(t) + c_2 y_2(t) + y_1(t) \sum \frac{\prod_{j=a}^{t-1} \beta(j) \sum \frac{r(t)y_1(t+1)} {\prod_{j=a}^{t}\beta(j)}} {y_1(t)y_1(t+1)}
\end{equation}
or
\begin{equation}\label{odiffform22}
y(t)=c_1y_1(t)+c_2y_2(t)+y_2(t)\sum\frac{\prod_{j=a}^{t-1}\beta(j)\sum\frac{r(t)y_2(t+1)}{\prod_{j=a}^{t}\beta(j)}}{y_2(t)y_2(t+1)},
\end{equation}
where $c_1$ and $c_2$ are arbitrary constants.
\end{theorem}

In \cite[Section 3]{jbb} the authors applied their alternative formula \eqref{jbbform} to give a general treatment of \eqref{ode} in the case of constant coefficient functions $p(t)\equiv 2P$ and $q(t)\equiv Q$ with $P$ and $Q$ real constants with $Q\neq 0$ and $P^2\neq Q$. We follow suit by giving a general treatment of \eqref{elaydiform} with constant coefficients using the results above from Theorem \ref{resultz2}. To wit, we analyze the ordinary difference equation
\begin{equation}\label{examplez2}
 y(t+2)+2\alpha y(t+1)+\beta y(t)=r(t), \quad t\in\Z,
\end{equation}
where $\alpha$ and $\beta$ are real constants with $\beta\neq 0$ and $\alpha^2\neq \beta$. The solution of \eqref{examplez2} can be written as
$$ y(t) = c_1 \left(-\alpha-\sqrt{\alpha^2-\beta}\right)^t + c_2 \left(-\alpha+\sqrt{\alpha^2-\beta}\right)^t+y_d. $$
From Theorem \ref{resultz2} we know that a particular solution to the inhomogeneous equation \eqref{examplez2} has the form
$$ y_d(t) = y_i(t) \sum_{x=a}^{t-1} \frac{\prod_{j=a}^{x-1}\beta\sum_{s=a}^{x-1}\frac{r(s)y_i(s+1)} {\prod_{j=a}^{s}\beta}}{y_i(x)y_i(x+1)}, \quad i=1,2, $$
which simplifies to
$$ y_d(t) = \sum_{x=a}^{t-1}\sum_{s=a}^{x-1} r(s)\left(-\alpha\pm\sqrt{\alpha^2-\beta}\right)^{t+s-2x}\beta^{x-1-s}.  $$
Thus, complete expressions for the solution of \eqref{examplez2} are given by
$$ y(t) = c_1 \lambda_1^t + c_2 \lambda_2^t + \sum_{x=a}^{t-1}\sum_{s=a}^{x-1} r(s)\lambda_1^{t+s-2x}\beta^{x-1-s} $$
or
$$ y(t) = c_1 \lambda_1^t + c_2 \lambda_2^t + \sum_{x=a}^{t-1}\sum_{s=a}^{x-1} r(s)\lambda_2^{t+s-2x}\beta^{x-1-s}, $$
where for simplicity we have taken $\lambda_1=-\alpha-\sqrt{\alpha^2-\beta}$ and $\lambda_2=-\alpha+\sqrt{\alpha^2-\beta}$.

\begin{remark}
A general treatment of \eqref{odyne} with constant coefficients on arbitrary time scales and even quantum equations is made difficult by the fact that, as seen in \eqref{ydform}, even with $p$ and $q$ taken to be constant functions, the term $-p+\mu q$ is not constant except for the very special cases of $\T=\R$ and $\T=h\Z$.
\end{remark}




\begin{thebibliography}{99}

\bibitem{ravi} R. P. Agarwal, 
\emph{Difference Equations and Inequalities: Theory, Methods, and Applications}, Second Edition, Revised and Expanded, Marcel Dekker, New York, 2000.

\bibitem{blest} D. Blest,
Reduction of order as a neglected method for the solution of inhomogeneous linear ordinary differential equations,
\emph{Int. J. Math. Educ. Sci. Technol.}, 21(3) (1990) 393--401.

\bibitem{bp1} M. Bohner and A. Peterson, 
\emph{Dynamic Equations on Time Scales, An Introduction with Applications}, Birkh\"auser, Boston, 2001.

\bibitem{bd} W. E. Boyce and R. C. DiPrima,
\emph{Elementary Differential Equations}, Seventh Edition, John Wiley \& Sons, New York, 2001.

\bibitem{coddington} E. A. Coddington,
\emph{An Introduction to Ordinary Differential Equations}, Prenctice--Hall Inc., Englewood Cliffs, 1961.

\bibitem{elaydi} S. N. Elaydi,
\emph{An Introduction to Difference Equations}, Springer, New York, 1996.

\bibitem{hilger} S. Hilger, 
Analysis on measure chains -- a unified approach to continuous and discrete calculus, \emph{Results Math.}, 18 (1990) 18--56.

\bibitem{hille} E. Hille,
\emph{Ordinary Differential Equations in the Complex Domain}, Dover, New York, 1997.

\bibitem{ince} E. L. Ince,
\emph{Ordinary Differential Equations}, Dover, New York, 1956.

\bibitem{jahnke} 
\emph{A History of Analysis}, translated from the German, edited by Hans Niels Jahnke, \emph{History of
Mathematics}, 24, American Mathematical Society, Providence, RI, London Mathematical Society, London, 2003.

\bibitem{jbb} P. Johnson, K. Busawon, and J. P. Barbot,
Alternative solution of the inhomogeneous linear differential equation of order two,
\emph{J. Math. Anal. Appl.}, 339 (2008) 582--589. 

\bibitem{kp2} W. Kelley and A. Peterson,
\emph{Difference Equations: An Introduction with Applications}, Second Edition, Academic Press, San Diego, 2000.

\bibitem{kp1} W. Kelley and A. Peterson,
\emph{The Theory of Differential Equations: Classical and Qualitative}, 
Pearson Prentice Hall, Upper Saddle River, NJ, 2004.

\bibitem{yosida} K. Yosida,
\emph{Lectures on Differential and Integral Equations}, Interscience Publishers, New York, 1960.

\end{thebibliography}
\end{document}